\documentclass[11pt]{amsart}
\usepackage{a4wide,latexsym,amssymb,amsthm,amsmath,amsfonts,url}

\begin{document}

\title{CR-structures of codimension $2$ on tangent bundles in Riemann-Finsler geometry}

\markboth{{\small\it {\hspace{5cm} CR-structures on tangent bundles in Riemann-Finsler geometry}}}{\small\it{CR-structures on tangent bundles in Riemann-Finsler geometry \hspace{5cm}}}

\author{Mircea Crasmareanu and Laurian-Ioan Pi\c scoran}
\date{April 11, 2013}

\maketitle

\begin{quote} {\small {\bf Abstract.}  We determine a $2$-codimensional CR-structure on the slit tangent bundle $T_0M$ of a Finsler manifold
$(M, F)$ by imposing a condition regarding the almost complex structure $\Psi $ associated to $F$ when restricted to the structural distribution of a framed $f$-structure. This condition is satisfied when $(M, F)$ is of scalar flag curvature (particularly flat) and in the Riemannian case $(M, g)$ this last condition means that $g$ is of constant curvature. This CR-structure is finally generalized by using one positive number but under more difficult conditions.
}
\end{quote}

\section*{Introduction}

The Finsler geometry is very rich in remarkable tensor fields $F$ of $(1, 1)$-type and associated structures. More precisely, there are: an (almost) tangent structure ($F^2=0$), an almost complex one ($F^2=-I$) and also an almost product structure $(F^2=I)$. In \cite{m:a} another well-known type of structures, namely $f$-structure ($F^3+F=0$) is obtained in this geometry. In fact, this $f$-structure belongs to a very interesting particular case which is called {\it framed $f$-structure} and has, in addition to $F$, a set of vector fields and differential $1$-forms interrelated. Moreover, a conformal deformation of the Sasaki type metric can be added in order to obtain a {\it metric framed $f$-structure}. This metric framed $f$-structure of M. Anastasiei was recently generalized in \cite{c:i} and \cite{p:z}.

\medskip

The present note is concerning with another kind of structures, namely the {\it CR-structures}, with an important r\^ole at the border between differential geometry and complex analysis, as it is pointed out in \cite{d:t}. We restrict ourselves at the real case; more precisely, based on a relationship between framed $f$-structures and CR-structure established in \cite[p. 130]{a:b} we found a CR-structure on the slit tangent bundle $T_0M$ of a Finsler manifold $(M, F)$. This CR-structure is constructed with the above almost complex structure denoted $\Psi _F$ in Section 2 and its existence is constrained by one condition expressing the vanishing of the Nijenhuis tensor of $\Psi _F$ on the structural distribution of the framed $f$-structure from \cite{m:a}. The above condition is expressed as a relation between the curvature of the Cartan nonlinear connection and the Jacobi endomorphism and is satisfied in dimension two or if $(M, F)$ is of scalar flag curvature which in the particular case of Riemannian geometry $(M, g)$ means that the metric $g$ has a constant curvature. Several important classes of Finsler manifolds with scalar flag curvature are discussed in Chapter 7 of \cite{c:z}.

\medskip

Inspired by \cite{p:z} we generalize this CR-structure using a real number $\beta >\frac{1}{2}$ but with more difficult conditions. More precisely, we take into account the same vector fields and $1$-forms as in the previous framed $f$-structure but deform the metric and the almost complex structure on both horizontal and vertical directions. At $\beta =1$ we recover the previous CR-structure.

\medskip

Finally, let us note that our CR-structures are of codimension $2$ and the (complex) geometry of these structures was studied in \cite{m:i1}-\cite{m:i2} and recently in \cite{k:z} and \cite{m:n}. But for the Riemannian case the only studies until now are on hypersurfaces of Sasakian manifolds (\cite{mu:i1}-\cite{mu:i2}) and not on (slit) tangent bundle.

\section{CR-structures from framed $f$-structures}

Framed $f$-structures constitute a particular case of $f$-structures and a detailed study of this class of tensor fields of $(1, 1)$-type, especially from a local point of view, can be found in \cite{y:k}.

\medskip

Let $N$ be a smooth $(2n+s)$-dimensional manifold with $n, s\geq 1$ and fix $D$ a distribution on $N$ of dimension $2n$. Considering $D$ as a vector bundle over $N$ let $\Gamma (D)$ be the module of its sections. Supposing $D$ is endowed with a morphism $J:D\rightarrow D$ of vector bundles satisfying $J^2=-I$ where $I$ is the identity (Kronecker) morphism on $D$, the pair $(D, J)$ is called {\it almost complex distribution}.

\medskip

The first main notion is given by \cite[p. 128]{a:b}:

\medskip

{\bf Definition 1.1} If for all $X, Y\in \Gamma (D)$ we have:
$$
\left\{
  \begin{array}{ll}
    [JX, JY]-[X, Y]\in \Gamma (D) \\
    N_J(X, Y):=[JX, JY]-[X, Y]-J([X, JY]+[JX, Y])=0
  \end{array}
\right. \eqno(1.1)
$$
then $(D, J)$ is a {\it CR-structure} on $N$ and the triple $(N, D, J)$ is a {\it CR-manifold}.

\smallskip

A second main notion is:

\medskip

{\bf Definition 1.2} Let $\varphi $ be a tensor field of $(1, 1)$-type and $s$ pairs $(\xi _a, \eta ^a)$, $1\leq a\leq s$ of vector fields and $1$-forms on $N$. If: \\
i) $\varphi ^3+\varphi =0$, $rank \varphi =2n$, \\
ii) $\varphi ^2=-I+\sum _{a=1}^s\eta ^a\otimes \xi _a$, $\varphi (\xi _a)=0$, $\eta ^a(\xi _b)=\delta ^a_b$, $\eta ^a\circ \varphi =0$, \\
then the data $(\varphi , \xi _a, \eta ^a)$ is called a {\it framed $f$-structure}.

\medskip

To a framed $f$-structure we associate \cite[p. 130]{a:b}:\\
1) {\it the torsion tensor field} $S$ of $(1, 2)$-type:
$$
S=N_{\varphi }+2\sum _{a=1}^sd\eta ^a\otimes \xi _a. \eqno(1.2)
$$
2) {\it the structural distribution} $D$:
$$
D=\{X\in \Gamma (TM); \eta ^1(X)=...=\eta ^s(X)=0\}=\cap _{a=1}^s\ker \eta ^a . \eqno(1.3)
$$
For a $1$-form $\eta $ we use the differential:
$$
2d\eta (X, Y)=X(\eta (Y))-Y(\eta (X))-\eta ([X, Y]). \eqno(1.4)
$$

\smallskip

These notions lead to:

\medskip

{\bf Definition 1.3} The framed $f$-structure is called $D$-{\it normal} if $S$ vanishes on $D$ i.e. $S(X, Y)=0$ for all $X, Y\in \Gamma (D)$.

\medskip

The relationship between the above structures was pointed out by A. Bejancu in Proposition 1.1 of \cite[p. 130]{a:b}:

\medskip

{\bf Proposition 1.4} {\it If $(\varphi , \xi _a, \eta ^a)$ is a $D$-normal framed $f$-structure then $(D, J=\varphi |_D)$ is a CR-structure}.

\medskip

{\bf Proof} The restriction $J$ of $\varphi $ to $D$ is obviously an almost complex structure. The conditions $(1.1)$ result from the fact that for $X, Y\in \Gamma (D)$ we have:
$$
S(X, Y)=0=[JX, JY]+\varphi ^2([X, Y])-\varphi ([X, JY]+[JX, Y])-\sum _{a=1}^s\eta ^a([X, Y])\xi _a. \eqno(1.5)
$$
For other details see the cited reference. \quad $\Box $

\section{A metric framed $f$-structure on the tangent bundle of a Finsler manifold}

Let $M$ be now a smooth $m$-dimensional manifold with $m\geq 2$ and $\pi : TM\rightarrow M$ its tangent bundle. Let $x=(x^i)=(x^1,..., x^m)$ be the local coordinates on $M$ and $(x, y)=(x^i, y^i)=(x^1,...,x^m, y^1,...., y^m)$ the induced local coordinates on $TM$. Denote by $O$ the null-section of $\pi $.

\medskip

Recall after \cite{c:z} that a {\it Finsler fundamental function} on $M$ is a map $F:TM\rightarrow \mathbb{R}_{+}$ with the following properties: \\
F1) $F$ is smooth on the slit tangent bundle $T_0M:=TM\setminus O$ and continuous on $O$, \\
F2) $F$ is positive homogeneous of degree $1$: $F(x,  \lambda y)=\lambda F(x, y)$ for every $\lambda >0$, \\
F3) the matrix $(g_{ij})=\left(\frac{1}{2}\frac{\partial ^2F^2}{\partial y^i\partial y^j}\right)$ is invertible and its associated quadratic form is positive definite. \\
The tensor field $g=\{g_{ij}(x, y); 1\leq i, j\leq \}$ is called {\it the Finsler metric} and the homogeneity of $F$ implies:
$$
F^2(x, y)=g_{ij}y^iy^j=y_iy^i \eqno(2.1)
$$
where $y_i=g_{ij}y^j$. The pair $(M, F)$ is called {\it Finsler manifold}.

\medskip

On $T_0M$ we have two distributions: \\
i) $V(TM):=\ker \pi _*$, called {\it the vertical distribution} and not depending of $F$. It is integrable and has the basis $\{\frac{\partial }{\partial y^i}; 1\leq i\leq m\}$. A remarkable section of it is {\it the Liouville vector field} $\Gamma =y^i\frac{\partial }{\partial y^i}$. \\
ii) $H(TM)$ with the basis $\{\frac{\delta }{\delta x^i}:=\frac{\partial }{\partial y^i}-N^j_i\frac{\partial }{\partial y^j}\}$ where:
$$
N^i_j=\frac{1}{2}\frac{\gamma ^i_{00}}{\partial y^j} \eqno(2.2)
$$
with $\gamma ^i_{00}=\gamma ^i_{jk}y^jy^k$ built from the usual Christoffel symbols:
$$
\gamma ^i_{jk}=\frac{1}{2}g^{ia}\left(\frac{\partial g_{ak}}{\partial x^j}+\frac{\partial g_{ja}}{\partial x^k}-\frac{\partial g_{jk}}{\partial x^a}\right). \eqno(2.3)
$$
$H(TM)$ is often called the {\it Cartan} (or canonical) {\it nonlinear connection} of the geometry $(M, F)$ and a remarkable section of it is {\it the geodesic spray}:
$$
S_F=y^i\frac{\delta }{\delta x^i}. \eqno(2.4)
$$
In particular, if $g$ does not depends on $y$ we recover the Riemannian geometry.

\medskip

The dual basis of the above local basis $\{\frac{\delta }{\delta x^i}, \frac{\partial }{\partial y^i}\}$ of $\Gamma (T_0M)$ is
$(dx^i, \delta y^i=dy^i+N^i_jdx^j)$. On $T_0M$ we have a Riemannian metric of Sasaki type:
$$
G_F=g_{ij}dx^i\otimes dx^j+g_{ij}\delta y^i\otimes \delta y^j. \eqno(2.5)
$$

\medskip

Another Finslerian object is the tensor field of $(1, 1)$-type $\Psi _F:\Gamma (T_0M)\rightarrow \Gamma (T_0M)$:
$$
\Psi _F\left(\frac{\delta }{\delta x^i}\right)=-\frac{\partial }{\partial y^i}, \quad \Psi _F\left(\frac{\partial }{\partial y^i}\right)=\frac{\delta }{\delta x^i}. \eqno(2.6)
$$
It results that $\Psi _F$ is an almost complex structure and the pair $(\Psi _F, G_F)$ is an almost K\"ahler structure on $T_0M$.

\medskip

In order to obtain a framed $f$-structure on $T_0M$ associated to the Finslerian function $F$, the following objects are considered in \cite{m:a}:
$$
\left\{
  \begin{array}{llll}
    \xi _1=S_F, \xi _2=\Gamma  \\
    \eta ^1=\frac{1}{F^2}y_idx^i, \eta ^2=\frac{1}{F^2}y_i\delta y^i \\
    \varphi =\Psi _F+\eta ^1\otimes \xi _2-\eta ^2\otimes \xi _1 \\
    G=\frac{1}{F^2}G_F.
  \end{array}
\right. \eqno(2.7)
$$
Then the main result of \cite{m:a} is that the data $(\varphi , \xi _1, \xi _2, \eta ^1, \eta ^2)$ is a framed $f$-structure on $T_0M$ with
$\eta ^a$ the $G$-dual of $\xi _a$, $1\leq a\leq 2$ and, moreover:
$$
G(\varphi \cdot, \varphi \cdot)=G-\eta ^1\otimes \eta ^1-\eta ^2\otimes \eta ^2. \eqno(2.8)
$$
Also, $\xi _a$ are unitary vector fields with respect to $G$ and $(G, \varphi , \xi _a, \eta ^a)$ is a {\it metric framed $f$-structure}.

\section{Putting all together}

The last paragraph of the previous Section provides the ingredients of the first Section with $N=T_0M$,  $s=2$ and $n=m-1$ which motivates our choice $m\geq 2$. The structural distribution is then:
$$
D_F=\ker \eta ^1\cap \ker \eta ^2=\{\xi _1\}^{\bot G}\cap \{\xi _2\}^{\bot G}=\{\xi _1\}^{\bot G_F}\cap \{\xi _2\}^{\bot G_F} \eqno(3.1)
$$
where $\{X\}^{\bot G}$ is the $G$-orthogonal complement of $span\{X\}$. We have $D_F=(span\{\xi _1, \xi_2\})^{\bot G_F}$ and this implies that $D_F$ has the dimension $2m-2$. For a geometrical meaning of the distribution $span\{\xi _1, \xi _2\}$ in \cite{m:a} is defined the differential $2$-form $\omega _F$, naturally associated to the metric framed $f$-structure:
$$
\omega _F=G(\cdot  , \varphi \cdot)  \eqno(3.2)
$$
and it follows that $span\{\xi _1, \xi _2\}$ is the kernel of $\omega _F$. Also, the homogeneity of $F$ implies the homogeneity of $S_F=\xi _1$ which means:
$$
[\Gamma , S_F]=[\xi _2, \xi _1]=\xi _1 \eqno(3.3)
$$
and thus $span\{\xi _1, \xi _2\}$ is an integrable distribution; see also Theorem 3.15 of \cite[p. 236]{b:f}.

\medskip

A concrete expression of $D_F$ appears in \cite[p. 11]{b:m}. More precisely, consider after the cited paper: \\
i) the horizontal vector fields:
$$
h_i=\frac{\delta }{\delta x^i}-\frac{1}{F^2}y_iS_F \eqno(3.4)
$$
and the corresponding $(m-1)$-distribution $\mathcal{H}_{m-1}=span\{h_i; 1\leq i\leq m\}$, \\
ii) the vertical vector fields:
$$
v_i=\frac{\partial }{\partial y^i}-\frac{1}{F^2}y_i\Gamma \eqno(3.5)
$$
and also the corresponding $(m-1)$-distribution $\mathcal{V}_{m-1}=span\{v_i; 1\leq i\leq m\}$. \\
We have:
$$
D_F=\mathcal{H}_{m-1}\otimes \mathcal{V}_{m-1} \eqno(3.6)
$$
and the same Theorem 3.15 of \cite[p. 236]{b:f} proves the integrability of $\mathcal{V}_{m-1}$; see also \cite[p. 12]{b:m}.

\medskip

Regarding the integrability of the nonlinear connection $H(TM)$ we have:
$$
\left[\frac{\delta }{\delta x^j}, \frac{\delta }{\delta x^k}\right]=R^i_{jk}\frac{\partial }{\partial y^i} \eqno(3.7)
$$
where:
$$
R^i_{jk}=\frac{\delta N^i_j}{\delta x^k}-\frac{\delta N^i_k}{\delta x^j}. \eqno(3.8)
$$

\medskip

The tensor field $R=\{R^i_{jk}(x, y); 1\leq i, j, k\leq m\}$ is called {\it the curvature} of the Cartan nonlinear connection and:
$$
R^i_j:=R^i_{kj}y^k \eqno(3.9)
$$
are the components of {\it the Jacobi endomorphism} $\Phi =R^i_{j}\frac{\partial }{\partial y^i}\otimes dx^j$, \cite[p. 5]{b:m}. We are ready now for the first main result:

\medskip

{\bf Theorem 3.1} {\it If the curvature tensor of $(M, F)$ has the form:
$$
R^i_{jk}=\lambda (X^i_ky_j-X^i_jy_k) \eqno(3.10)
$$
with $\lambda $ a smooth function on $T_0M$ and the tensor field $\{X^i_j(x, y); 1\leq i, j\leq m\}$ satisfying:
$$
y_iX^i_j=y_j \eqno(3.11)
$$
for all $i, j\in \{1,...,m\}$ then the pair $(D_F, J_F=\Psi _F|_{D_F})$ is a CR-structure on $T_0M$ }.

\medskip

{\bf Proof} We express the Nijenhuis tensor field of $\Psi _F$ as:
$$
N_{\Psi _F}(X, Y)=[\Psi _FX, \Psi _FY]-[X, Y]-\Psi _F(A(X, Y))=B(X, Y)-\Psi _F(A(X, Y)) \eqno(3.12)
$$
with $A(X, Y):=[X, \Psi _FY]+[\Psi _FX, Y]$ and $B(X, Y)=[\Psi _FX, \Psi _FY]-[X, Y]$. It follows that $B(X, Y)=A(\Psi _FX, Y)$ and then:
$$
N_{\Psi _F}(X, Y)=A(\Psi _FX, Y)-\Psi _F\circ A(X, Y). \eqno(3.13)
$$

We prove firstly that $A$ is a $D_F$-valued $(0, 2)$-tensor field. From $(3.7)$ and:
$$
\left[\frac{\delta }{\delta x^j}, \frac{\partial }{\partial y^k}\right]=\frac{\partial N^i_j}{\partial y^k}\frac{\partial }{\partial y^i}=
\frac{\partial ^2 \gamma ^i_{00}}{\partial y^j\partial y^k}\frac{\partial }{\partial y^i} \eqno(3.14)
$$
we obtain:
$$
A\left(\frac{\delta }{\delta x^j}, \frac{\delta }{\delta x^k}\right)=A\left(\frac{\partial }{\partial y^j}, \frac{\partial }{\partial y^k}\right)=0, \quad A\left(\frac{\delta }{\delta x^j}, \frac{\partial }{\partial y^k}\right)=R^i_{jk}\frac{\partial }{\partial y^i} \eqno(3.15)
$$
which means that $\eta ^1\circ A=0$ and:
$$
A=R^i_{jk}dx^j\wedge \delta y^k\otimes \frac{\partial }{\partial y^i}. \eqno(3.16)
$$
A main identity in Finsler geometry is:
$$
y_iR^i_{ab}=0 \eqno(3.17)
$$
and then $\eta ^2\circ A=0$ which conclude the first part of the proof.

\medskip

Secondly, we search for the framework of Proposition 1.4. The torsion tensor $S$ on $D_F$ is:
$$
S(X, Y)=N_{\varphi }(X, Y)-\eta ^1([X, Y])\xi _1-\eta ^2([X, Y])\xi _2
$$
with:
$$
N_{\varphi }(X, Y)=[\Psi _FX, \Psi _FY]+\varphi ^2([X, Y])-\varphi \circ A(X, Y).
$$
Since $\varphi $ is an element of a framed $f$-structure we get:
$$
N_{\varphi }(X, Y)=[\Psi _FX, \Psi _FY]-[X, Y]+\eta ^1([X, Y])\xi _1+\eta ^2([X, Y])\xi _2-\varphi \circ A(X, Y)
$$
and from the definition $(2.7_3)$ of $\varphi $ it follows:
$$
S(X, Y)=[\Psi _FX, \Psi _FY]-[X, Y]-(\Psi _F+\eta ^1\otimes \xi _2-\eta ^2\otimes \xi _1)\circ A(X, Y)=N_{\Psi _F}(X, Y). \eqno(3.18)
$$
In local coordinates we have:
$$
N_{\Psi _F}=R^i_{jk}\delta y^j\wedge \delta y^k\otimes \frac{\partial }{\partial y^i} \eqno(3.19)
$$
and then $N_{\Psi _F}$ has components only when applied on the pair $(v_a, v_b)$. A long but straightforward computation yields:
$$
N_{\Psi _F}(v_a, v_b)=2\left[R^i_{ab}+\frac{1}{F^2}(R^i_ay_b-R^i_by_a)\right]\frac{\partial }{\partial y^i} \eqno(3.20)
$$
and therefore the normality condition is:
$$
F^2R^i_{ab}=R^i_by_a-R^i_ay_b    \eqno(3.21)
$$
which can be expressed as:
$$
N_{\Psi _F}=\eta ^2\wedge \left(R^i_k\delta y^k\otimes \frac{\partial }{\partial y^i}\right). \eqno(3.22)
$$
The relation $(3.10)$ yields:
$$
R^i_k=\lambda (F^2X^i_k-y^aX^i_ay_k) \eqno(3.23)
$$
and then, both sides of $(3.21)$ are equal with $\lambda F^2(X^i_ky_j-X^i_jy_k)$ which gives the final conclusion. The condition $(3.11)$ corresponds to the relation $(3.17)$.

Let us also point out that the condition $(3.10)$ gives the following expression for the Nijenhuis tensor:
$$
N_{\Psi _F}=2\lambda F^2\eta ^2\wedge \left(X^i_j\delta y^j\otimes \frac{\partial }{\partial y^i}\right) \eqno(3.24)
$$
which yields again the vanishing of $N_{\Psi _F}$ on $D_F$ due to the presence of $\eta ^2$. Concerning the tensor field $A$ we have:
$$
A=\lambda F^2\left[\eta ^1\wedge (X^i_j\delta y^j\otimes \frac{\partial }{\partial y^i})-(X^i_jdx^j\otimes \frac{\partial }{\partial y^i})\wedge \eta ^2\right] \eqno(3.25)
$$
which proves the relations: $\eta ^1\circ A=\eta ^2\circ A=0$. \quad $\Box $

\medskip

{\bf Example 3.2} Recall that in dimension $2$ the Nijenhuis tensor field of any almost complex structure vanishes. Then every $2$-dimensional Finsler manifold $(M_2, F)$ satisfies the condition of Theorem 3.1. Let $V(TM)$ be spanned by the vector fields $\Gamma $ and $V$ respectively $H(TM)$ be spanned by the vector fields $S_F$ and $H$. Then $D_F$ is spanned by $V$ and $H$ and:
$$
J_F(H)=-V, \quad J_F(V)=H. \eqno(3.26)
$$
We have that $H$ is a linear combination of $h_1$ and $h_2$ while $V$ is a linear combination of $v_1$ and $v_2$. $\Box $

\medskip

In order to consider examples in any dimension we remark that a solution of condition $(3.11)$ is:
$$
X^i_j=\mu \delta ^i_j+(1-\mu )\frac{y^iy_j}{F^2} \eqno(3.27)
$$
again with $\mu $ a smooth function on $T_0M$. It follows:

\medskip

{\bf Example 3.3} If $\mu =1$ then $X^i_j=\delta ^i_j$ and the Finsler manifold $(M, F)$ is {\it of scalar flag curvature} $\lambda $ since:
$$
R^i_{jk}=\lambda (\delta ^i_ky_j-\delta ^i_jy_k) \eqno(3.28)
$$
and then:
$$
R^i_k=\lambda (\delta ^i_kF^2-y^iy_k). \eqno(3.29)
$$

\smallskip

{\bf Corollary 3.4} {\it If $(M, F)$ is of scalar flag curvature then $(D_F=(span\{S_F, \Gamma \})^{\bot G_F}, J_F)$ is a CR-structure on $T_0M$}.

\medskip

Remark also that the hypothesis of scalar flag curvature yields:
$$
N_{\Psi _F}=2\lambda F^2\eta ^2\wedge \pi _{V(TM)} \eqno(3.30)
$$
where $\pi _{V(TM)}$ is the projector on the vertical part in the $G_F$-orthogonal decomposition $T(T_0M)=H(TM)\oplus V(TM)$ i.e $\pi _{V(TM)}=\delta y^i\otimes \frac{\partial }{\partial y^i}$. However, $\Psi _F$ is integrable only in the flat case (i.e. $\lambda =0$) since $N_{\Psi _F}(\Gamma , v_a)=2\lambda F^2v_a$. The integrability of $\Psi _F$ as a tensor field of $(1, 1)$-type which is equivalent with the integrability of the Cartan nonlinear connection of $(M, F)$ and then $(T_0M, \Psi _F, G_F)$ is a K\"ahler manifold.

\medskip

{\bf Particular case 3.5} (Riemannian geometry) Let $g=(g_{ij}(x))$ be a Riemannian metric on $M$. Then $\gamma ^i_{jk}(x, y)=\Gamma ^i_{jk}(x)$ the Riemannian Christoffel symbols and:
$$
R ^i_{jk}(x, y)=R^i_{jka}(x)y^a \eqno(3.31)
$$
with $R_g=(R^i_{jka})$ the Riemannian curvature tensor of $g$. It results that a Riemannian geometry $(M, F=(g_{ij}(x)y^iy^j)^{\frac{1}{2}})$ is of scalar flag curvature if and only if $g$ is of constant curvature. Therefore on the slit tangent bundle of a space form $(M, g)$ there exists a CR-structure on the distribution complementary (with respect to the Sasaki lift of $g$) to the distribution generated by the Liouville vector field and the geodesic spray $S_g$. \quad $\Box $

\medskip

{\bf Example 3.6} Returning to the general non-Riemannian case $(3.27)$ with $\mu =0$ we get:
$$
X^i_j=\frac{y^iy_j}{F^2} \eqno(3.32)
$$
and then $R^i_{jk}=0$ which means that $(M, F)$ is flat, a situation belonging also to the Example 3.3 for vanishing scalar curvature. \quad $\Box $

For the general $\mu $ we have:
$$
N_{\Psi _F}=2\lambda F^2\eta ^2\wedge [\mu \pi _{V(TM)}+(1-\mu )\eta ^2\otimes \Gamma ] =2\lambda \mu F^2\eta ^2\wedge \mu \pi _{V(TM)}. \eqno(3.33)
$$

\section{A $1$-parametric generalization}

Let $\alpha >0$ and $\beta >0$ two positive numbers as well as the smooth function $v:[0, +\infty )\rightarrow \mathbb{R}$ which, following the approach of \cite{p:z}, will be considered as $v=v(\tau )$ with $\tau =F^2$. Supposing that:
$$
\alpha +2\tau v(\tau )>0 \eqno(4.1)
$$
for any $\tau $ in the cited paper is constructed the smooth function:
$$
w=-\frac{\beta v}{\alpha +\tau v} \eqno(4.2)
$$
and the Riemannian metric on $T_0M$:
$$
\bar{G}=G_{ij}dx^i\otimes dx^j+H_{ij}\delta y^i\otimes \delta y^j \eqno(4.3)
$$
where:
$$
\left\{
  \begin{array}{ll}
    G_{ij}=\frac{1}{\beta }g_{ij}+\frac{v}{\alpha \beta }y_iy_j \\
    H_{ij}=\beta g_{ij}+wy_iy_j.
  \end{array}
\right. \eqno(4.4)
$$
Inspired by \cite{p:z} we define also:
$$
\left\{
  \begin{array}{ll}
    \bar{\xi }_1=(\beta +w\tau )S_F, \quad \bar{\xi }_2=\Gamma =\xi _2  \\
    \bar{\eta }^1=\frac{1}{\tau }y_idx^i=\eta ^1, \quad \bar{\eta }^2=(\frac{\beta }{\tau }+w)y_i\delta y^i  \\ [2mm]
    \bar{\Psi }_F(\frac{\delta }{\delta x^i})=-G^a_i\frac{\partial }{\partial y^a}, \quad \bar{\Psi }_F(\frac{\partial }{\partial y^i})=H^a_i\frac{\delta }{\delta x^a}
  \end{array}
\right. \eqno(4.5)
$$
where the lift of indices in the third line is constructed with $g^{-1}=(g^{ab})$. In fact, the only difference between us and \cite{p:z} is with respect to $1$-form $\bar{\eta }^i$; in order to reobtain that of Section 2 we divide with $\tau $ the $1$-forms of Peyghan-Zhong. With a computation similar to that of Theorem 4.8 of Peyghan-Zhong we derive that $(\bar{G}, \bar{\varphi }, \bar{\xi }_a, \bar{\eta }^a)$ with:
$$
\bar{\varphi }=\bar{\Psi }_F+\bar{\eta }^1\otimes \bar{\xi }_2-\bar{\eta }^2\otimes \bar{\xi }_1 \eqno(4.6)
$$
is a metric framed $f$-structure on $T_0M$ if and only if:
$$
\beta +\tau w=1. \eqno(4.7)
$$
From this condition we get that $\bar{\xi }_a=\xi _a$ and $\bar{\eta }^a=\eta ^a$. From $(4.2)$ and $(4.7)$ we obtain:
$$
v(\tau )=\frac{\alpha (\beta -1)}{\tau }, \quad w(\tau ) =\frac{1-\beta}{\tau }. \eqno(4.8)
$$
In the particular case $\alpha =\beta =1$ we recover the metric framed $f$-structure of Anastasiei since $v=w\equiv 0$.

\medskip

Now, under condition $(4.7)$ we have the same structural distribution $D_F$ but the expression of the tensor field:
$$
\bar{A}(X, Y):=[X, \bar{\Psi }_FY]+[\bar{\Psi }_FX, Y] \eqno(4.9)
$$
is more complicated. More detailed:
$$
\left\{
  \begin{array}{ll}
    \bar{A}(\frac{\delta }{\delta x^j}, \frac{\delta }{\delta x^k})=\left(\frac{\delta G^v_j}{\delta x^k}-\frac{\delta G^v_k}{\delta x^j}+G^u_j\frac{\partial N^v_k}{\partial y^u}-G^u_k\frac{\partial N^v_j}{\partial y^u}\right)\frac{\partial }{\partial y^v} \\
    \bar{A}(\frac{\partial }{\partial y^j}, \frac{\partial }{\partial y^k})=\left(\frac{\partial H^v_k}{\partial y^j}-\frac{\partial H^v_j}{\partial y^k}\right)\frac{\delta }{\delta x^v}+\left(H^u_j\frac{\partial N^v_u}{\partial y^k}-H^u_k\frac{\partial N^v_u}{\partial y^j}\right)\frac{\partial }{\partial y^v} \\
    \bar{A}(\frac{\delta }{\delta x^j}, \frac{\partial }{\partial y^k})=\frac{\delta H^v_k}{\delta x^j}\frac{\delta }{\delta x^v}+\left(H^u_kR^v_{ju}+\frac{\partial G^v_j}{\partial y^k}\right)\frac{\partial }{\partial y^v}
  \end{array}
\right. \eqno(4.10)
$$
where, with $(4.7)$:
$$
\left\{
  \begin{array}{ll}
    G_{ij}=\frac{1}{\beta }g_{ij}+\frac{\beta -1}{\beta \tau }y_iy_j, \quad H_{ij}=\beta g_{ij}+\frac{1-\beta }{\tau }y_iy_j \\
    G^a_j=\frac{1}{\beta }\delta ^a_j+\frac{\beta -1}{\beta \tau }y^ay_j, \quad H^a_j=\beta \delta ^a_j+\frac{1-\beta }{\tau }y^ay_j \\
    \bar{\Psi }_F(\frac{\delta }{\delta x^i})=-\frac{1}{\beta }\frac{\partial }{\partial y^i}+\frac{1-\beta }{\beta \tau }y_i\Gamma , \quad \bar{\Psi }_F(\frac{\partial }{\partial y^i})=\beta \frac{\delta }{\delta x^i}+\frac{1-\beta }{\tau }y_iS_F.
  \end{array}
\right. \eqno(4.11)
$$
It results that $\alpha $ disappears and this motives the title of this Section, namely $1$-parametric ge\-ne\-ra\-li\-za\-tion and not $2$-parametric. Note that $\bar{\Psi }_F\left(h_i\right)=-\frac{1}{\beta }v_i$ and $\bar{\Psi }_F(v_i)=\beta h_i$.

\smallskip

Then:
$$
\left\{
  \begin{array}{ll}
    \bar{A}(\frac{\delta }{\delta x^j}, \frac{\delta }{\delta x^k})=\frac{\beta -1}{\beta \tau }\left[\frac{\delta }{\delta x^k}\left(y_jy^v\right)-\frac{\delta }{\delta x^j}\left(y_ky^v\right)\right]\frac{\partial }{\partial y^v} \\
    \bar{A}(\frac{\partial }{\partial y^j}, \frac{\partial }{\partial y^k})=(1-\beta )\left[\frac{\partial }{\partial y^j}(\frac{y_ky^v}{\tau })-\frac{\partial }{\partial y^k}(\frac{y_jy^v}{\tau })\right]\frac{\delta }{\delta x^v} \\
   \bar{A}(\frac{\delta }{\delta x^j}, \frac{\partial }{\partial y^k})=\frac{1-\beta }{\tau }\frac{\delta }{\delta x^j}\left(y_ky^v\right)\frac{\delta }{\delta x^v}+\left[\beta R^v_{jk}+\frac{1-\beta }{\tau }y_ky^uR^v_{ju}+\frac{\beta -1}{\beta }\frac{\partial }{\partial y^k}\left(\frac{y_jy^v}{\tau }\right)\right]\frac{\partial }{\partial y^v}.
  \end{array}
\right. \eqno(4.12)
$$

Choosing $\alpha =1$ the second main result is:

\medskip

{\bf Theorem 4.1} {\it Let $\beta >\frac{1}{2}$ and the smooth functions $v(\tau )=-w(\tau )=\frac{\beta -1}{\tau }$. If for any $X, Y\in D_F$ we have}:\\
1) $\bar{A}(X, Y)\in D_F$, \\
2) $N_{\bar{\Psi }_F}(X, Y)=0$,   \\
{\it then $(D_F, \bar{J}_F=\bar{\Psi }_F|_{D_F})$ is a CR-structure on} $T_0M$.

\medskip

{\bf Proof} The condition in $\beta $ is the expression of $(4.1)$. Exactly as in the proof of Theorem 3.1 we have:
$$
S(X, Y)=N_{\bar{\Psi }_F}(X, Y)-\eta ^1(\bar{A}(X, Y))\xi _2+\eta ^2(\bar{A}(X,  Y))\xi _1. \eqno(4.13)
$$
and the conclusion follows directly. Let us note that 1) corresponds to the condition $(1.1_1)$ while 2) corresponds to the condition $(1.1_2)$.  \quad $\Box $

\medskip

Let us remark that:
$$
\beta \eta ^2\circ \bar{A}\left(\frac{\delta }{\delta x^j}, \frac{\delta }{\delta x^k}\right)=\eta ^1\circ \bar{A}\left(\frac{\delta }{\delta x^j}, \frac{\partial }{\partial y^k}\right)-\eta ^1\circ \bar{A}\left(\frac{\delta }{\delta x^k}, \frac{\partial }{\partial y^j}\right). \eqno(4.14)
$$
and then, the vanishing of $\eta ^1\circ \bar{A}\left(\frac{\delta }{\delta x^a}, \frac{\partial }{\partial y^b}\right)$ implies the vanishing of $\eta ^2\circ \bar{A}\left(\frac{\delta }{\delta x^u}, \frac{\delta }{\delta x^v}\right)$. The vanishing of the former expression means that $y_k$ are eigenvalue for $\frac{\delta }{\delta x^j}$:
$$
\frac{\delta y_k}{\delta x^j}=\left(-\frac{N^a_jy_a}{F^2}\right)y_k \eqno(4.15)
$$
and then $y_k$ are eigenvalues for the geodesic spray:
$$
S_F(y_k)=\left(-\frac{N^a_jy^jy_a}{F^2}\right)y_k. \eqno(4.16)
$$
Such condition holds in the Euclidian space $(\mathbb{R}^m, g_{ij}=\delta _{ij})$ but here the expression $\eta ^2\circ \bar{A}(\frac{\delta }{\delta x^j}, \frac{\partial }{\partial y^k})$ is non-vanishing since:
$$
y_v\frac{\partial }{\partial y^k}\left(\frac{y_jy^v}{F^2}\right)=\delta _{jk}-\frac{y_jy^k}{F^2}\neq 0 \eqno(4.17)
$$
and then it remains an open problem to find Riemannian and/or Finsler manifolds satisfying the Theorem 4.1 with $\beta \neq 1$.

\medskip

{\bf Acknowledgement.} The first author has been supported by the research grant PN-II-ID-PCE-2012-4-0131.

\small{

\bigskip

\noindent Mircea Crasmareanu \hfill Laurian-Ioan Pi\c scoran \\
Faculty of Mathematics \hfill Department of Mathematics and Computer Science \\
University "Al. I. Cuza", 700506, Ia\c si, Rom\^ania \hfill North Univ. Center of Baia Mare, 430122\\
mcrasm@uaic.ro \hfill Technical Univ. of Cluj Napoca, Rom\^ania \\
http://www.math.uaic.ro/$\sim$mcrasm \hfill plaurian@yahoo.com
}


\begin{thebibliography}{Nov}

\bibitem{m:a} M. Anastasiei, {\it A framed $f$-structure on tangent manifold of a Finsler space}, An. Univ. Bucure\c sti Mat. Inform., 49(2000), no. 2, 3-9. MR1892254 (2003b:53082)

\bibitem{a:b} A. Bejancu, {\it Geometry of CR-submanifolds}, Mathematics and its Applications (East European Series), vol. 23, D. Reidel Publishing Co., Dordrecht, 1986. MR0861408 (87k:53126)

\bibitem{b:f} A. Bejancu; H. R. Farran, {\it Foliations and geometric structures}, Mathematics and Its Applications, vol. 580, Springer, Dordrecht, 2006. MR2190039 (2006j:53034)

\bibitem{b:m} I. Bucataru; Z. Muzsnay, {\it Projective and Finsler metrizability: parameterization-rigidity of the geodesics}, Internat. J. Math., 23(2012), no. 9, 1250099, 15 pp. MR2959445

\bibitem{c:z} S.-S. Chern; Z. Shen, {\it Riemann-Finsler geometry}, Nankai Tracts in Mathematics, 6. World Scientific Publishing Co. Pte. Ltd., Hackensack, NJ, 2005. MR2169595 (2006d:53094)

\bibitem{d:t} S. Dragomir; G. Tomassini, {\it Differential geometry and analysis on CR manifolds}, Progress in Mathematics, vol. 246. Birkhäuser Boston, Inc., Boston, MA, 2006. MR2214654 (2007b:32056)

\bibitem{c:i} C. Ida, {\it Some framed $f$-structures on transversally Finsler foliations}, Ann. Univ. Mariae Curie-Sklodowska Sect. A, 65(2011), no. 1, 87-96. MR2825153 (2012f:53157)

\bibitem{k:z} S.-Y. Kim; D. Zaitsev, {\it Equivalence and embedding problems for CR-structures of any codimension.}, Topology, 44(2005), no. 3, 557-584. MR2122216 (2005j:32037)

\bibitem{m:n} C. Medori; M. Nacinovich, {\it Standard CR manifolds of codimension} 2, Transform. Groups, 6(2001), no. 1, 53-78. MR1825168 (2002c:32055)

\bibitem{m:i1} R. I. Mizner, {\it CR structures of codimension $2$}, J. Differential Geom., 30(1989), no. 1, 167-190. MR1001274 (90h:32046)

\bibitem{m:i2} R. I. Mizner, {\it Almost CR structures, $f$-structures, almost product structures and associated connections}, Rocky Mountain J. Math., 23(1993), no. 4, 1337-1359. MR1256452 (95d:32016)

\bibitem{mu:i1} M.-I. Munteanu, {\it CR-structures of CR-codimension $2$ on hypersurfaces in Sasakian manifolds}, in "Differential geometry and its applications", 157-163, Matfyzpress, Prague, 2005. MR2268930 (2007h:53086)

\bibitem{mu:i2} M.-I. Munteanu, {\it New aspects on CR-structures of codimension $2$ on hypersurfaces of Sasakian manifolds}, Arch. Math. (Brno), 42(2006), no. 1, 69-84. MR2227114 (2007h:53066)

\bibitem{p:z} E. Peyghan; C. Zhong, {\it A framed $f$-structure on the tangent bundle of a Finsler manifold}, Ann. Polon. Math., 104(2012), no. 1, 23-41. MR2885972

\bibitem{y:k} K. Yano; M. Kon, {\it Structures on manifolds}, Series in Pure Mathematics, vol. 3, World Scientific Publishing Co., Singapore, 1984. MR0794310 (86g:53001)

\end{thebibliography}
\end{document}